\documentclass[fleqn,11pt]{article}
\usepackage{amsmath}
\usepackage{amsfonts}
\usepackage{graphicx}

\title{\bf Relatively normal-slant helices lying on a surface and their characterizations}
\author{ Nesibe MAC\.IT\thanks{ Department of Mathematics, Faculty of Science and Arts, Y{\i}ld{\i}z Technical University, \.Istanbul, Turkey,  \ Email: {\tt ngurhan@yildiz.edu.tr }}, \,\,
Mustafa D\"ULD\"UL\thanks{ Department of Mathematics, Faculty of Science and Arts, Y{\i}ld{\i}z Technical University, \.Istanbul, Turkey,  \ Email: {\tt mduldul@yildiz.edu.tr }}}
\date{}
\newtheorem{definition}{Definition}
\newtheorem{theorem}{Theorem}

\newtheorem{remark}{Remark}
\newtheorem{corollary}{Corollary}
\newtheorem{example}{Example}

\numberwithin{equation}{section}

\usepackage{footnote}
\usepackage{graphics}
\usepackage{multicol}
\usepackage{booktabs}
\usepackage[flushleft]{threeparttable}
\usepackage[font=small,labelfont=bf]{caption}

\begin{document}
\maketitle

\begin{abstract}
In this paper, we consider a regular curve on an oriented surface in Euclidean 3-space with the Darboux frame $\{\mathsf{T},\mathsf{V},\mathsf{U}\}$ along the curve, where $\mathsf{T}$ is the unit tangent vector field of the curve, $\mathsf{U}$ is the surface normal restricted to the curve and $\mathsf{V}=\mathsf{U}\times \mathsf{T}$. We define a new curve on a surface by using the Darboux frame. This new curve whose vector field $\mathsf{V}$ makes a constant angle with a fixed direction is called as relatively normal-slant helix. We give some characterizations for such curves and obtain their axis. Besides we give some relations between some special curves (general helices, integral curves, etc.) and relatively normal-slant helices. Moreover, when a regular surface is given by its implicit or parametric equation, we introduce the method for generating the relatively normal-slant helix with the chosen direction and constant angle on the given surface.
\end{abstract}

\section{Introduction}

Helical curves play important roles in not only CAD and CAGD but also in science and nature. In the field of computer aided design and computer graphics, helices can be
used for the tool path description, the simulation of kinematic motion or the design
of highways, etc. [24]. They also arise in the structure of DNA, fractal geometry (e.g. hyperhelices), etc. [5, 23]. In recent years, helical curves are studied widely in Euclidean spaces (e.g. [1-3, 11, 15, 20]) and in non-Euclidean spaces (e.g. [4, 8, 12]).

A space curve whose tangent vector makes a constant angle with a fixed direction is called a generalized helix and characterized by Lancret's theorem which says the ratio of torsion to curvature is constant [21]. Generalized helices are also studied in higher dimensional spaces [19, 22].

In 2004, Izumiya and Takeuchi define a slant helix in $\mathbb{E}^3$ by the property that the principal normal vector makes a constant angle with a fixed direction. These curves are characterized by the constancy of geodesic curvature function of the principal normal indicatrix of the curve [10]. In 2005 Kula and Yayl{\i} have studied spherical indicatrices of a slant helix and obtained that the spherical images are spherical helix [14].

Besides, slant helices according to Bishop frame are studied by Bukcu and Karacan in 2009 [6]. Also, slant helices according to quaternionic frame [13] and slant helices in 3-dimensional Lie group [17] are investigated in 2013.

In [25], Z{\i}plar and Yayl{\i} introduce a Darboux helix which is defined as a curve whose Darboux vector $\omega=\tau \mathsf{T} + \kappa \mathsf{B}$ makes a constant angle with a fixed direction, and give some characterizations of such curves in 2012. 

Puig-Pey et al. introduce the method for generating the general helix for both implicit and parametric surfaces [20].

In a recent paper, Do\u{g}an and Yayl{\i} study isophote curves and their characterizations in Euclidean 3-space [7]. An isophote curve is defined as a curve on a surface whose unit normal vector field restricted to the curve makes a constant angle with a fixed direction. They also obtain the axis of an isophote curve.

In this study, we define a relatively normal-slant helix on a surface by using the Darboux frame $\{\mathsf{T},\mathsf{V},\mathsf{U}\}$ along the curve whose vector field $\mathsf{V}$ makes a constant angle with a fixed direction. We give some characterizations for such curves and obtain their axis. Besides we give some relations between some special curves (general helices, integral curves, etc.) and relatively normal-slant helices. Moreover, when a regular surface is given by its implicit or parametric equation, we introduce the method for generating the relatively normal-slant helix with the chosen direction and constant angle on the given surface. This method is based on an initial-value problem of 1st order ordinary differential equations.

This paper is organized as follows: Section 2 includes some basic definitions. We define relatively normal-slant helices in section 3 and give some characterizations for such curves. The method for finding the relatively normal-slant helix on a given surface is introduced in section 4 for parametric and implicit surfaces and its application is presented in section 5. Section 6 includes some relations between the special curves and relatively normal-slant helices. Some further characterizations related with the subject are also given in section 7.

\section{Preliminaries}\setcounter{equation}{0}
Let $M$ be an oriented surface, and $\alpha: I\subset\mathbb{R}\rightarrow M$ be a regular curve with arc-length parametrization. If the Frenet frame along the curve is denoted by $\{\mathsf{T}, \mathsf{N}, \mathsf{B}\}$, the Frenet formulas are given by
\[
\mathsf{T}'=\kappa \mathsf{N},\quad \mathsf{N}'=-\kappa \mathsf{T}+\tau \mathsf{B}, \quad \mathsf{B}'=-\tau \mathsf{N},
\]
where $\mathsf{T}$ is unit tangent vector, $\mathsf{N}$ is principal normal vector, $\mathsf{B}$ is the binormal vector; $\kappa$ and $\tau$ are the curvature and the torsion of $\alpha$, respectively. On the other hand, if we denote the Darboux frame along the curve $\alpha$ by $\{\mathsf{T}, \mathsf{V}, \mathsf{U}\}$, we have the derivative formulae of the Darboux frame as
\[
\mathsf{T}'=\kappa_g \mathsf{V}+\kappa_n\mathsf{U},\quad \mathsf{V}'=-\kappa_g \mathsf{T}+\tau_g \mathsf{U}, \quad \mathsf{U}'=-\kappa_n\mathsf{T}-\tau_g \mathsf{V},
\]
where $\mathsf{T}$ is the unit tangent vector of the curve, $\mathsf{U}$ is the unit normal vector of the surface restricted to the curve, $\mathsf{V}$ is the unit vector given by $\mathsf{V}=\mathsf{U}\times \mathsf{T}$, and $\kappa_g$, $\kappa_n$, $\tau_g$ denote the geodesic curvature, normal curvature, geodesic torsion of the curve, respectively [18].

The relations between geodesic curvature, normal curvature, geodesic torsion
and $\kappa$, $\tau$ are given as follows:
\[
\kappa_g=\kappa \sin\varphi,\qquad \kappa_n=\kappa \cos\varphi,\qquad \tau_g=\tau+\frac{d\varphi}{ds},
\]
where $\varphi$ is the angle between the vectors $\mathsf{N}$ and $\mathsf{U}$.

If the surface $M$ is given by its parametric equation
\begin{equation*}
X(u,v)=\big(x(u,v),y(u,v),z(u,v)\big),
\end{equation*}
we may write $\alpha(s)=X(u(s),v(s))$ for the curve $\alpha$. Thus, the tangent vector of the curve $\alpha$ at the point $\alpha(s)$ becomes
\begin{equation}
\frac{d\alpha}{ds}=X_u\frac{du}{ds}+X_v\frac{dv}{ds},
\end{equation}
where $X_u$ and $X_v$ denote the partial differentiation of $X$ with respect to $u$ and $v$, respectively. Hence, since $\alpha(s)$ is a unit-speed curve, we have
\begin{equation}
E\left(\frac{du}{ds}\right)^2+2F\frac{du}{ds}\frac{dv}{ds}+G\left(\frac{dv}{ds}\right)^2=1,
\end{equation}
where $E=\langle X_u,X_u\rangle$, $F=\langle X_u,X_v\rangle$ and $G=\langle X_v,X_v\rangle$ are the first fundamental form coefficients of the surface.

On the other hand, if the surface $M$ is given by its implicit equation $f(x,y,z)=0$, then the unit speed curve $\alpha(s)=\big(x(s),y(s),z(s)\big)$ lying on $M$ satisfies
\begin{equation}
f_x\frac{dx}{ds}+f_y\frac{dy}{ds}+f_z\frac{dz}{ds}=0
\end{equation}
and
\begin{equation}
\left(\frac{dx}{ds}\right)^2+\left(\frac{dy}{ds}\right)^2+\left(\frac{dz}{ds}\right)^2=1,
\end{equation}
where $f_a=\frac{\partial f}{\partial a}$.

\begin{definition}[Slant helix]
A unit speed curve is called a slant helix if its unit principal normal vector makes a constant angle with a fixed direction [10].
\end{definition}
\begin{theorem}
Let $\alpha$ be a unit speed curve with $\kappa\not=0$. Then $\alpha$ is a slant helix if and only if
\begin{equation*}
\sigma(s)=\left(\frac{\kappa^2}{\left(\kappa^2+\tau^2\right)^{3/2}}\left(\frac{\tau}{\kappa}\right)'\right)(s)
\end{equation*}
is a constant function [10].
\end{theorem}

\begin{definition}[$\mathsf{V}$-direction curve]
Let $\gamma$ be a unit speed curve on an oriented surface $M$ and $\{\mathsf{T},\mathsf{V},\mathsf{U}\}$ be the Darboux frame along $\gamma$. The curve $\overline{\gamma}$ lying on $M$ is called a $\mathsf{V}$-direction curve of $\gamma$ if it is the integral curve of the vector field $\mathsf{V}$. In other words, if  $\overline{\gamma}$ is the $\mathsf{V}$-direction curve of $\gamma$, then
$\mathsf{V}(s)=\overline{\gamma}'(s)$ [16].
\end{definition}
\begin{definition}[Darboux vector fields]
Let $M$ be an oriented surface, and $\alpha: I\subset\mathbb{R}\rightarrow M$ be a regular curve with arc-length parametrization. The vector fields $D_n(s), D_r(s), D_o(s)$ along $\alpha$ given by
\begin{equation*}
D_n(s)=-\kappa_n(s)\mathsf{V}(s)+\kappa_g(s)\mathsf{U}(s),
\end{equation*}
\begin{equation*}
D_r(s)=\tau_g(s)\mathsf{T}(s)+\kappa_g(s)\mathsf{U}(s),
\end{equation*}
\begin{equation*}
D_o(s)=\tau_g(s)\mathsf{T}(s)-\kappa_n(s)\mathsf{V}(s)
\end{equation*}
are called the normal Darboux vector field, the rectifying Darboux vector field and the osculating Darboux vector field along $\alpha$, respectively [9].
\end{definition}

\section{Relatively normal-slant helix and its axis}

In this section, we introduce a new kind of helix which is lying on a surface.

Let $M$ be an oriented surface in $\mathbb{E}^3$ and $\gamma$ be a regular curve lying on $M$. Let us denote the Darboux frame along $\gamma$ with $\{\mathsf{T},\mathsf{V},\mathsf{U}\},$ where $\mathsf{T}$ is the unit tangent vector field of $\gamma,$ $\mathsf{U}$ is the unit normal vector field of the surface which is restricted to the curve $\gamma$ and $\mathsf{V}=\mathsf{U}\times \mathsf{T}$.

\begin{definition}[Relatively normal-indicatrix]
Let $\gamma$ be a unit speed curve (with arc-length parameter $s$) on an oriented surface $M$ and $\{\mathsf{T},\mathsf{V},\mathsf{U}\}$ be the Darboux frame along $\gamma$. We define the curve which $\mathsf{V}(s)$ draws on the unit sphere $\mathbb{S}^2$ as the relatively normal-indicatrix of $\gamma$, and denote it by $\gamma_v$. Thus, $\gamma_v(s)=\mathsf{V}(s)$.
\end{definition}

\begin{definition}[Relatively normal-slant helix]
Let $\gamma$ be a unit speed curve on an oriented surface $M$ and $\{\mathsf{T},\mathsf{V},\mathsf{U}\}$ be the Darboux frame along $\gamma$. The curve $\gamma$ is called a relatively normal-slant helix if the vector field $\mathsf{V}$ of $\gamma$ makes a constant angle with a fixed direction, i.e. there exists a fixed unit vector $d$ and a constant angle $\theta$ such that $\langle \mathsf{V},d\rangle=cos\theta.$
\end{definition}

\begin{theorem}
A unit speed curve $\gamma$ on a surface $M$ with $(\tau_g(s),\kappa_g(s))\neq (0,0)$ is a relatively normal-slant helix if and only if
\begin{equation}
\sigma_v(s)=\left(\frac{1}{\left(\kappa_g^2 +\tau_g^2\right)^{\frac{3}{2}}}\Big(\tau_g'\kappa_g-\kappa_g' \tau_g- \kappa_n\big(\kappa_g^2+ \tau_g^2\big)\Big)\right)(s)
\end{equation}
is a constant function.
\end{theorem}
{\bf Proof.}
$(\Rightarrow)$ Let the unit speed curve $\gamma$ on a surface $M$ with $(\tau_g(s),\kappa_g(s))\neq (0,0)$ be a relatively normal-slant helix. Then, by the definition, the vector field $\mathsf{V}$ of $\gamma$ makes a constant angle with a fixed direction. Thus, the relatively normal-indicatrix of the curve $\gamma$ is part of a circle, i.e.
it has constant curvature and zero torsion.

For the relatively normal-indicatrix of the curve $\gamma$, we may write $\gamma_v(s)=\mathsf{V}(s)$ which, by differentiating with respect to $s$, yields
\[
\gamma_v'(s)=-\kappa_g\mathsf{T}+\tau_g\mathsf{U}
\]
\[
\gamma_v''(s)=\big(-\kappa_g'-\kappa_n\tau_g\big)\mathsf{T}-\big(\kappa_g^2+ \tau_g^2\big)\mathsf{V}+\big(-\kappa_g\kappa_n+\tau_g'\big)\mathsf{U}.
\]
The curvature $\kappa_v$ and the torsion $\tau_v$ of $\gamma_v$ are obtained as
\begin{equation}
\left\{\begin{array}{l}
\kappa_v(s)=\frac{\parallel \gamma_v'\times\gamma_v''\parallel}{\parallel\gamma_v'\parallel^3}
=\sqrt{1+\big(\sigma_v(s)\big)^2},\\
\\
\tau_v(s)=\frac{1}{\left(\kappa_g^2(s) +\tau_g^2(s)\right)^{\frac{11}{2}}}\frac{\sigma_v'(s)}{(1+(\sigma_v(s))^{2})},
\end{array}\right.
\end{equation}
where
$\sigma_v(s)=\frac{1}{\left(\kappa_g^2 +\tau_g^2\right)^{\frac{3}{2}}}\Big(\tau_g'\kappa_g-\kappa_g' \tau_g- \kappa_n\big(\kappa_g^2+ \tau_g^2\big)\Big)$. Hence, since $\gamma_v(s)$ has constant curvature and zero torsion, we have $\sigma_v(s)=constant$.

$(\Leftarrow)$ Suppose that $\sigma_v(s)=constant$. Then, by using (3.2) the relatively normal-indicatrix of $\gamma$ has constant curvature and zero torsion, i.e. relatively normal-indicatrix is part of a circle which means $\mathsf{V}$ makes constant angle with a fixed direction. Thus, $\gamma$ is a relatively normal-slant helix.

\begin{remark}
The constant function $\sigma_v(s)$ is also the geodesic curvature of the relatively normal-indicatrix of $\gamma$.
\end{remark}
\begin{remark}
The invariant $\sigma_v(s)$ is equal to $\frac{-\delta_r(s)}{\sqrt{\kappa_g^2 +\tau_g^2}}$, where $\delta_r$ is given in [9].
\end{remark}

We may give the following corollaries for the special cases of the relatively normal-slant helices.
\begin{corollary}
{\bf i)}
Let $\gamma$ be an asymptotic curve on $M$ with $\kappa_g(s)\not=0$. Then, $\gamma$ is a relatively normal-slant helix on $M$ if and only if
\begin{equation}
\left(\frac{\kappa_g^2}{\left(\kappa_g^2 +\tau_g^2\right)^{\frac{3}{2}}}\left(\frac{\tau_g}{\kappa_g}\right)'\right)(s)=constant.
\end{equation}
{\bf ii)} If $\gamma$ is an asymptotic curve on $M$ with $\kappa_g(s)\not=0$, we have $\kappa=\kappa_g, \tau=\tau_g$. Then, $\gamma$ is a relatively normal-slant helix on $M$ if and only if $\gamma$ is a slant helix.
\end{corollary}

\begin{corollary}
Let $\gamma$ be a geodesic curve on $M$ with $\tau_g(s)\not=0$. Then, $\gamma$ is a relatively normal-slant helix on $M$ if and only if
\begin{equation}
\frac{\kappa_n(s)}{\tau_g(s)}=constant,
\end{equation}
i.e. $\gamma$ is a relatively normal-slant helix on $M$ if and only if $\gamma$ is a generalized helix on $M$.
\end{corollary}

\begin{corollary}
Let $\gamma$ be a line of curvature on $M$ with $\kappa_g(s)\not=0$. Then, $\gamma$ is a relatively normal-slant helix on $M$ if and only if
\begin{equation}
\frac{\kappa_n(s)}{\kappa_g(s)}=constant.
\end{equation}
\end{corollary}
\vspace{0.7cm}

Now, it is usual to ask that how can we find the axis of a relatively normal-slant helix. Our goal is now to obtain the axis of our new defined slant helix.

Let a unit speed curve $\gamma$ lying on an oriented surface $M$ be a relatively normal-slant helix. Then, by the definition, $\mathsf{V}$ makes a constant angle $\theta$ with the fixed unit vector $d=a\mathsf{T}+b\mathsf{V}+c\mathsf{U}$, i.e.
$\langle \mathsf{V},d\rangle=\cos\theta=b$. By differentiating the last equation with respect to $s$ gives
\begin{equation*}
\langle \mathsf{V}',d\rangle=0\quad\Rightarrow\quad \langle -\kappa_g \mathsf{T}+\tau_g \mathsf{U},d\rangle=0.
\end{equation*}
If we differentiate $d=a\mathsf{T}+b\mathsf{V}+c\mathsf{U}$ with respect to $s$, we obtain
\begin{equation*}
d'=(a'-b\kappa_g-c\kappa_n)\mathsf{T}+(a\kappa_g-c\tau_g)\mathsf{V}+(c'+a\kappa_n+b\tau_g)\mathsf{U}=0,
\end{equation*}
i.e.
\begin{equation}
\left\{
\begin{array}{l}
a'-b\kappa_g-c\kappa_n=0,\\
a\kappa_g-c\tau_g=0,\\
c'+a\kappa_n+b\tau_g=0.
\end{array}\right.
\end{equation}
If we substitute $c=\frac{\kappa_g}{\tau_g}a, \tau_g(s)\not=0,$ into the first and third equations and multiply the third one with $\frac{\kappa_g}{\tau_g}$, we obtain a first order homogeneous linear differential equation as
\[
\left(1+\left(\frac{\kappa_g}{\tau_g}\right)^2\right)a'+\frac{\kappa_g}{\tau_g}\left(\frac{\kappa_g}{\tau_g}\right)'a=0
\]
which has the general solution
\[
a=\lambda\frac{\tau_g}{\sqrt{\kappa_g^2+\tau_g^2}},
\]
where $\lambda$ is the constant of integration.
Substituting this solution into the second equation of (3.6), we obtain
\[
c=\lambda\frac{\kappa_g}{\sqrt{\kappa_g^2+\tau_g^2}}
\]
On the other hand, since $d$ is unit length, i.e. $a^2+b^2+c^2=1$,
we find $\lambda=\pm\sin\theta$.
Therefore, $d$ can be written as
\begin{equation}
d=\pm\frac{\tau_g}{\sqrt{\kappa_g^2+\tau_g^2}}\sin\theta \mathsf{T}+\cos\theta \mathsf{V}\pm\frac{\kappa_g}{\sqrt{\kappa_g^2+\tau_g^2}}\sin\theta \mathsf{U}.
\end{equation}

We need to determine the constant angle $\theta$ to complete the axis.

If we differentiate $\langle \mathsf{V}',d\rangle=0$ with respect to s along the curve $\gamma$, we obtain
\[
\langle \mathsf{V}'',d\rangle=\pm\left(\frac{\tau_g'\kappa_g-\kappa_g' \tau_g- \kappa_n(\kappa_g^2+ \tau_g^2)}{\sqrt{\kappa_g^2 +\tau_g^2}}\right)\sin\theta-(\kappa_g^2 +\tau_g^2)\cos\theta=0.
\]
As a result we have
\begin{equation}
\cot \theta=\pm\frac{\tau_g'\kappa_g-\kappa_g' \tau_g- \kappa_n(\kappa_g^2+ \tau_g^2)}{(\kappa_g^2 +\tau_g^2)^\frac{3}{2}}.
\end{equation}
Hence, finding $\theta$ from (3.8) and substituting the result into (3.7) gives us the axis of the relatively normal-slant helix (By using (3.8), it is easy to see that $d$ is a constant vector, i.e. $d'=0$).

\section{Calculating a relatively normal-slant helix on a surface}
In this section we introduce the methods for finding the relatively normal-slant helix on a given surface. We discuss the methods separately for parametric and implicit surfaces.

\subsection{Relatively normal-slant helix on a parametric surface}

Let $M$ be a regular oriented surface in $\mathbb{E}^3$ with the parametrization $X=X(u,v)$. Our goal is now, when a fixed unit direction $d$ and a constant angle $\theta$ are given, to give the method which enables us to find the relatively normal-slant helix (if exists) lying on $M$ which accepts $d$ as an axis and $\theta$ as the constant angle.

Let $\gamma(s)=X(u(s),v(s))$ be the unit speed relatively normal-slant helix lying on $M$ with axis $d$, constant angle $\theta$, and Darboux frame field $\{\mathsf{T},\mathsf{V},\mathsf{U}\}$. Thus, we have $\langle \mathsf{V},d\rangle= \cos\theta$. We need to find $u(s), v(s)$ to obtain $\gamma$.

Since
\begin{equation*}
\gamma'=\mathsf{T}=X_u\frac{du}{ds}+X_v\frac{dv}{ds}\quad\text{and}\quad
\mathsf{U}=\frac{X_u\times X_v}{\parallel X_u\times X_v\parallel},
\end{equation*}
we have
\[
\mathsf{V}=\mathsf{U}\times \mathsf{T}=\frac{X_u\times X_v}{\parallel X_u\times X_v\parallel}\times \left(X_u\frac{du}{ds}+X_v\frac{dv}{ds}\right)
\]
i.e.
\begin{equation*}
\mathsf{V}=\frac{1}{\parallel X_u\times X_v\parallel}\left[\left(EX_v-FX_u\right)\frac{du}{ds}+\left(FX_v-GX_u\right)\frac{dv}{ds} \right].
\end{equation*}
Substituting the last equation into $\langle \mathsf{V},d\rangle= \cos\theta$ gives us
\begin{equation}
\left(E\langle X_v,d\rangle-F\langle X_u,d\rangle\right)\frac{du}{ds}+\left(F\langle X_v,d\rangle-G\langle X_u,d\rangle\right)\frac{dv}{ds}= \cos\theta\parallel X_u\times X_v\parallel.
\end{equation}
Combining Eqs. (2.2) and (4.1), we obtain
\begin{equation}
\left\{
\begin{array}{l l }
{\frac{du}{ds}=\frac{2\cos\theta(EG-F^2)^{\frac{3}{2}}\langle X_v,d\rangle\pm\sqrt{\Delta}}{2A(EG-F^2)}}\\
\\
{\frac{dv}{ds}=\frac{-2\cos\theta(EG-F^2)^{\frac{3}{2}}\langle X_u,d\rangle\mp\sqrt{\Delta^*}}{2A(EG-F^2)}}
\end{array}
\right.
\end{equation}
where $A$, $\Delta$, and $\Delta^*$ are given by
\[
A=E\langle X_v,d\rangle^2-2F\langle X_u,d\rangle\langle X_v,d\rangle+G\langle X_u,d\rangle^2,
\]
\[
\Delta=4\cos^2\theta(EG-F^2)^2\left[\langle X_v,d\rangle^2(EG-F^2)-AG\right]
\]
\[
\qquad+4A(EG-F^2)\left[F\langle X_v,d\rangle-G\langle X_u,d\rangle\right]^2,
\]
\[
\Delta^*=4\cos^2\theta(EG-F^2)^2\left[\langle X_u,d\rangle^2(EG-F^2)-AE\right]
\]
\[
\qquad+4A(EG-F^2)\left[E\langle X_v,d\rangle-F\langle X_u,d\rangle\right]^2.
\]
If we solve the system (4.2) together with the initial point
\begin{equation}
\left\{
\begin{array}{l l }
u(0)=u_0\\
v(0)=v_0,
\end{array}
\right.
\end{equation}
we obtain the desired relatively normal-slant helix on $M$ by substituting $u(s), v(s)$ into $X(u,v)$.

\begin{remark}
{\bf i)} If $\Delta$ and/or $\Delta^*$ is negative, it means there does not exist a relatively normal-slant helix with the given axis and angle.\\
{\bf ii)} If $\Delta\geq 0$ and $\Delta^*\geq 0$, then we have two relatively normal-slant helices on the surface.
\end{remark}

\subsection{Relatively normal-slant helix on an implicit surface}
Let $M$ be a surface given in implicit form by $f(x,y,z)=0$. Let us now find the relatively normal-slant helix $\gamma(s)$ which makes the given constant angle $\theta$ with the given axis $d=(a,b,c)$ and lying on $M$.

Let $\gamma(s)=(x(s),y(s),z(s))$ and $\{\mathsf{T},\mathsf{V},\mathsf{U}\}$ be its Darboux frame field. We need to find $x(s), y(s), z(s)$ to obtain $\gamma(s)$ (We assume that $s$ is the arc-length parameter).

Since the unit normal vector field of the surface is
$\mathsf{U}=\frac{\nabla f}{\parallel\nabla f\parallel}$, we obtain
\[
\mathsf{V}=\frac{\nabla f}{\parallel\nabla f\parallel}\times \mathsf{T}
=\frac{1}{\parallel\nabla f\parallel}\left( f_y\frac{dz}{ds}-f_z\frac{dy}{ds},f_z\frac{dx}{ds}-f_x\frac{dz}{ds},f_x\frac{dy}{ds}-f_y\frac{dx}{ds} \right).
\]
If we substitute the last equation into $\langle \mathsf{V},d\rangle= \cos\theta$, we get
\begin{equation}
(bf_z-cf_y)\frac{dx}{ds}+(cf_x-af_z)\frac{dy}{ds}+(af_y-bf_x)\frac{dz}{ds}=\parallel\nabla f\parallel\cos\theta.
\end{equation}
If we consider (4.4) and (2.3), we obtain $\frac{dx}{ds}, \frac{dy}{ds}$ by means of $\frac{dz}{ds}$ as
\begin{equation}
\frac{dx}{ds}=\frac{1}{\Omega}\left[f_y(af_y-bf_x)-f_z(cf_x-af_z))\frac{dz}{ds}-f_y\parallel\nabla f\parallel\cos\theta \right],
\end{equation}
\begin{equation}
\frac{dy}{ds}=\frac{1}{\Omega}\left[f_z(bf_z-cf_y)-f_x(af_y-bf_x))\frac{dz}{ds}+f_x\parallel\nabla f\parallel\cos\theta \right],
\end{equation}
where $\Omega=cf_x^2-af_xf_z-bf_yf_z+cf_y^2\neq 0$.
Substituting (4.5) and (4.6) into (2.4) give us a quadratic equation with respect to $\frac{dz}{ds}$ as
\[
q_1\left(\frac{dz}{ds}\right)^2+q_2\frac{dz}{ds}+q_3=0,
\]
where
\[
q_1=\frac{1}{\Omega^2}\Big[b^2f_x^4+a^2f_y^4+(a^2+b^2)f_z^4
-2abf_xf_y^3-2acf_xf_z^3-2bcf_yf_z^3
\]
\[
\qquad-2abf_yf_x^3+(a^2+b^2)f_x^2f_y^2 
+(c^2+2b^2)f_x^2f_z^2+(2a^2+c^2)f_y^2f_z^2
\]
\[
\qquad-4abf_xf_yf_z^2\Big]+1,
\]
\[
q_2=\frac{1}{\Omega^2}\Big[2\parallel\nabla f\parallel\cos\theta(bf_xf_z^2-af_yf_x^2+bf_xf_y^2-af_yf_z^2+bf_x^3-af_y^3)\Big],
\]
\[
q_3=\frac{1}{\Omega^2}\Big[\parallel\nabla f\parallel^2\cos^2\theta(f_x^2+f_y^2)\Big]-1.
\]
From this equation we have
\begin{equation}
\frac{dz}{ds}=\frac{-q_2\pm \sqrt{q_2^2-4q_1q_3}}{2q_1}.
\end{equation}
If we substitute (4.7) into (4.5) and (4.6), we obtain an explicit 1st order ordinary differential equation system. Thus, together with the initial point
\begin{equation}
\left\{
\begin{array}{l l }
x(0)=x_0\\
y(0)=y_0\\
z(0)=z_0
\end{array}
\right.
\end{equation}
we have an initial value problem. The solution of this problem gives us the relatively normal-slant helix on $M$.
\begin{remark}
{\bf i)} If $q_2^2-4q_1q_3<0$ at the point $(x_0,y_0,z_0)$, then there does not exist any relatively normal-slant helix with the given direction $d$ and angle $\theta$.\\
{\bf ii)} If $q_2^2-4q_1q_3>0$ at the point $(x_0,y_0,z_0)$, then we have two relatively normal-slant helices passing through the initial point.
\end{remark}

\section{Examples}

\begin{example}
Let $M$ be the surface given by $X(u,v)=(u\cos v, u\sin v, u^2)$. Applying our method, the relatively normal-slant helix on $M$ which makes the constant angle $\theta=\frac{\pi}{3}$ with the chosen direction $d=(0,0,1)$ and starting from the inital point $P=(1,0,1)$ is given in Figure 1. Another application is given in Figure 2 in which we obtain two relatively normal-slant helices lying on the surface $(x^2+y^2)z^2+\frac{x^2+y^2}{4}-\frac{1}{4}=0$ with the initial point $\left(\frac{-1}{\sqrt{13}},0,-\sqrt 3\right)$
(in each figure the general helices are obtained by the method given in [20]).
\begin{figure}[h]
\centering
  \includegraphics[width=4in]{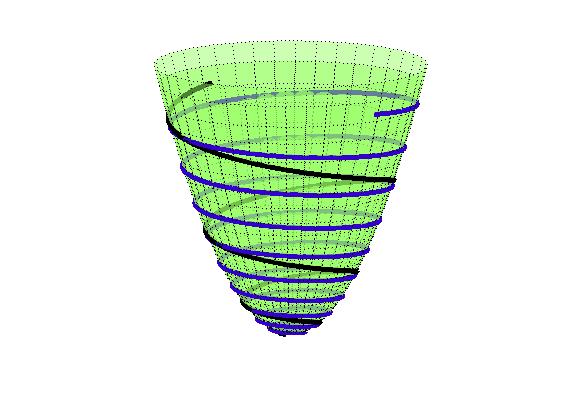}
  \caption{General helix (colored with blue) and relatively normal-slant helix (colored with black) with same axis $d=(0,0,1)$ and angle $\theta=\pi/3$}
\end{figure}
\begin{figure}[h]
\centering
  \includegraphics[width=4in]{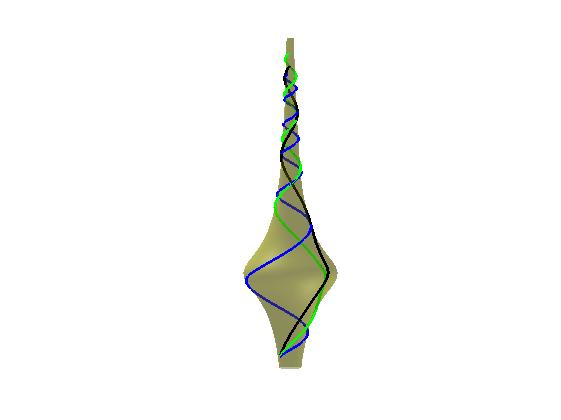}
  \caption{General helix (colored with blue ($\theta=\pi/3$)) and relatively normal-slant helices (colored with black ($\theta=\pi/3$) and green ($\theta=\pi/4$)) with same axis $d=(0,0,1)$ }
\end{figure}
\end{example}

\section{Some special curves related with a relatively normal-slant helix}

\begin{theorem}
Let $M$ be an oriented surface in $\mathbb{E}^3$, $\gamma$ be a unit speed curve on $M$, and $\overline{\gamma}$ be the $\mathsf{V}$-direction curve of $\gamma$. Then $\gamma$ is a relatively normal-slant helix if and only if $\overline{\gamma}$ is a general helix.
\end{theorem}
{\bf Proof.}
Since $\overline{\gamma}$ is the $\mathsf{V}$-direction curve of $\gamma$,
the relations between the curvature $(\overline{\kappa})$ and the torsion $(\overline{\tau})$ of $\overline{\gamma}$ and the geodesic curvature, normal curvature and geodesic torsion of $\gamma$ are given by, [16],
\begin{equation}
\overline{\kappa}=\sqrt{\kappa_g^2+\tau_g^2},\quad
\overline{\tau}=-\kappa_n+\frac{\kappa_g \tau_g'-\kappa_g'\tau_g}{\kappa_g^2+\tau_g^2}.
\end{equation}
Thus, we have
\begin{equation*}
\frac{\overline{\tau}}{\overline{\kappa}}=\left(\frac{1}{(\kappa_g^2 +\tau_g^2)^{\frac{3}{2}}}(\tau_g'\kappa_g-\kappa_g' \tau_g- \kappa_n(\kappa_g^2+ \tau_g^2)\right).
\end{equation*}
It means that $\gamma$ is a relatively normal-slant helix if and only if $\overline{\gamma}$ is a general helix.

\begin{theorem}
Let $M$ be an oriented surface in $\mathbb{E}^3$, $\gamma: I\subset\mathbb{R}\rightarrow M$ be a regular curve with arc-length parametrization and $D_r$ be the rectifying Darboux vector field of $\gamma$. Then $\gamma$ is a relatively normal-slant helix if and only if the integral curve of the rectifying Darboux vector field $D_r$ is a circular helix.
\end{theorem}
{\bf Proof.}
The rectifying Darboux vector field $D_r$ of $\gamma$ is given by
\begin{equation}
D_r(s)=\tau_g(s)\mathsf{T}(s)+\kappa_g(s)\mathsf{U}(s).
\end{equation}
Let $\beta$ be the integral curve of $D_r$. So we have
\[
D_r=\beta'=\tau_g\mathsf{T}+\kappa_g\mathsf{U}.
\]
By differentiating this equation with respect to $s$ we obtain
\[
\beta''=(\tau_g'-\kappa_g\kappa_n)\mathsf{T}+(\kappa_g'+\kappa_n\tau_g)\mathsf{U},
\]
\[
\beta'''=(\tau_g'-\kappa_g\kappa_n)'\mathsf{T}+(\tau_g'-\kappa_g\kappa_n)(\kappa_g\mathsf{V}+\kappa_n\mathsf{U})
\]
\[
\qquad+(\kappa_g'+\kappa_n\tau_g)'\mathsf{U}+(\kappa_g'+\kappa_n\tau_g)(-\kappa_n\mathsf{T}-\tau_g\mathsf{V}).
\]
Therefore, since
\[
\beta'\times\beta''=\Big(-\kappa_n(\kappa_g^2+\tau_g^2)+\tau_g'\kappa_g
-\tau_g\kappa_g'\Big)\mathsf{V},
\]
the curvature $\kappa_\beta$ of $\beta$ is obtained as
\[
\kappa_\beta=\frac{\parallel\beta'\times\beta''\parallel}{\parallel\beta'\parallel^3}=\left(\frac{1}{(\kappa_g^2 +\tau_g^2)^{\frac{3}{2}}}(\tau_g'\kappa_g-\kappa_g' \tau_g- \kappa_n(\kappa_g^2+ \tau_g^2)\right)=\sigma_v(s)
\]
and the torsion $\tau_\beta$ of $\beta$ is obtained as
\[
\tau_\beta=\frac{\langle\beta'\times\beta'',\beta'''\rangle}{\parallel\beta'\times\beta''\parallel^2}=1.
\]
Therefore
\[
\frac{\tau_\beta}{\kappa_\beta}=\frac{1}{\sigma_v(s)}.
\]
It means that $\gamma$ is a relatively normal-slant helix if and only if the integral curve of the rectifying Darboux vector field $D_r$ is a circular helix.

\begin{theorem}
Let $\gamma: I\rightarrow M\subset\mathbb{E}^3$ be a unit speed curve with $(\kappa_g(s),\tau_g(s))\not=(0,0)$. Then the following conditions are equivalent.\medskip

i) $\gamma$ is a relatively normal-slant helix,\medskip

ii) The relatively normal-indicatrix of $\gamma$ is part of a circle in $\mathbb{S}^2$,\medskip

iii) $\sigma_v(s)=\left(\frac{1}{(\kappa_g^2 +\tau_g^2)^{\frac{3}{2}}}(\tau_g'\kappa_g-\kappa_g' \tau_g- \kappa_n(\kappa_g^2+ \tau_g^2)\right)(s)$ is constant,\medskip

iv) There exists a unit vector $d$ such that $\left\langle d,\frac{D_r}{||D_r||}\right\rangle$ is constant,\medskip

v) The $\mathsf{V}$-direction curve of $\gamma$ is a general helix,\medskip

vi) The $D_r$-integral curve of $\gamma$ is a circular helix.

\end{theorem}

\section{Further characterizations}
Similar to the previous sections, we can also give the following characterizations for the general helices and isophote curves.

\begin{theorem}
Let $M$ be an oriented surface in $\mathbb{E}^3$, $\gamma: I\subset\mathbb{R}\rightarrow M$ be a regular curve with arc-length parametrization and $D_n$ be the normal Darboux vector field of $\gamma$. Then $\gamma$ is a general helix if and only if the integral curve of the normal Darboux vector field $D_n$ is a circular helix.
\end{theorem}
\begin{theorem}
Let $M$ be an oriented surface in $\mathbb{E}^3$, $\gamma: I\subset\mathbb{R}\rightarrow M$ be a regular curve with arc-length parametrization and $D_o$ be the osculating Darboux vector field of $\gamma$. Then $\gamma$ is an isophote curve if and only if the integral curve of the osculating Darboux vector field $D_o$ is a circular helix.
\end{theorem}
{\bf Acknowledgment.} The authors would like to thank the reviewer for the valuable comments and suggestions.

\end{document}